\documentclass{article}

\usepackage{amssymb,amsmath,amsthm}

\newtheorem{theorem}{Theorem}

\newtheorem{proposition}[theorem]{Proposition}
\newtheorem{remark}[theorem]{Remark}

\newenvironment{keywords}{%
\begin{center}
\begin{minipage}[c]{12cm}%
 {\bf Keywords:}}%
 {\end{minipage}
\end{center}}

\newenvironment{msc}{%
\begin{center}%
\begin{minipage}[c]{12cm}%
{\bf 2000 Mathematics Subject Classification:}}%
{\end{minipage}
\end{center}
\medskip}

\begin{document}

\title{Two-dimensional Newton's Problem\\ of Minimal Resistance\footnote{Presented
at the 4th Junior European Meeting on ``Control and
Optimization'', Institute of Mathematics and Physics, Bia\l ystok
Technical University, Bia\l ystok, Poland, 11-14 September 2005.
Research report CM06/I-01. Accepted (06-July-2006) to Control \&
Cybernetics.}}

\author{Cristiana J. Silva\\
\texttt{cjoaosilva@mat.ua.pt} \and
Delfim F. M. Torres\\
\texttt{delfim@mat.ua.pt}}
\date{Department of Mathematics\\
University of Aveiro\\
3810-193 Aveiro, Portugal}

\maketitle

%%%%%%%%%%%%%%%%%%%%%%%%%%%%%%%%%%%%%%%%%%%%%%%%%

\begin{abstract}
Newton's problem of minimal resistance is one of the first
problems of optimal control: it was proposed, and its solution
given, by Isaac Newton in his masterful \emph{Principia
Mathematica}, in 1686. The problem consists of determining, in
dimension three, the shape of an axis-symmetric body, with
assigned radius and height, which offers minimum resistance when
it is moving in a resistant medium. The problem has a very rich
history and is well documented in the literature. Of course, at
first glance, one suspects that the two dimensional case should be
well known. Nevertheless, we have looked into numerous references
and ask at least as many experts on the problem, and we have not
been able to identify a single source. Solution was always
plausible to everyone who thought about the problem, and writing
it down was always thought not to be worthwhile. Here we show that
this is not the case: the two-dimensional problem is more rich
than the classical one, being, in some sense, more interesting.
Novelties include: (i) while in the classical three-dimensional
problem only the restricted case makes sense (without restriction
on the monotonicity of admissible functions the problem doesn't
admit a local minimum), we prove that in dimension two the
unrestricted problem is also well-posed when the ratio height
versus radius of base is greater than a given quantity; (ii) while
in three dimensions the (restricted) problem has a unique
solution, we show that in the restricted two-dimensional problem
the minimizer is not always unique -- when the height of the body
is less or equal than its base radius, there exists infinitely
many minimizing functions.
\end{abstract}


\begin{msc}
49K05 (76G25, 76M30).
\end{msc}

%%%%%%%%%%%%%%%%%%%%%%%%%%%%%%%%%%%%%%%%%

\section{Introduction}

Newton's aerodynamical problem, in dimension three, is a classic
problem (see \textrm{e.g.}
\cite{CD:Amaral:1913,CD:Veubeke:1966,CD:Kneser:1913}). It consists
in joining two given points of the plane by a curve's arc that,
while turning around a given axis, generate the body of revolution
offering the least resistance when moving in a fluid in the
direction of the axis. Newton has considered several hypotheses:
that the body moves with constant velocity, and without rotation,
on a very rare and homogeneous medium of particles which are all
equal; that the axis-symmetric body is inscribed in a cylinder of
height $H$ and radius $r$; that the particles of the medium are
infinitesimally small and immovable (there exists no temperature
motion of particles); that collisions of the particles with the
body are absolutely elastic. Newton has indicated in the
\emph{Mathematical principles of natural philosophy} the correct
solution to his problem. He has not explained, however: how such
solution can be obtained; how the problem is formulated in the
language of mathematics. This has been the work of many
mathematicians since Newton's time (see \textrm{e.g.}
\cite{MR56:4953,CD:Tikhomirov:1990,CD:arXiv:math.OC/0404237}).
Extensions of Newton's problem is a topic of current intensive
research, with many questions remaining open challenging problems.
Recent results, obtained by relaxing Newton's hypotheses, include:
non-symmetric bodies \cite{CD:BK:1993}; one-collision non-convex
bodies \cite{CD:CL:2001}; collisions with friction
\cite{CD:friction:2002}; multiple collisions allowed
\cite{CD:P1:2003}; temperature noise of particles
\cite{CD:arXiv:math.OC/0404194,math.OC/0407406}. Here we are
interested in the classical problem, under the classical
hypotheses considered by Newton. Our main objective is to study
the apparently more simpler Newton's problem of minimal resistance
for a two-dimensional body moving with constant velocity in a
homogeneous rarefied medium of particles. The first work on a
two-dimensional Newton-type problem seems to be
\cite{CD:arXiv:math.OC/0404194}, where the authors study the
problem in a chaotically moving media of particles (in the
classical problem particles are immovable). The results in
\cite{CD:arXiv:math.OC/0404194} were later generalized to
dimension three \cite{math.OC/0407406}. This paper is motivated by
the results in \cite{math.OC/0407406}: when one considers
temperature motion of particles, the three-dimensional problem
admits only two types of solutions; while the two-dimensional case
is more rich, showing solutions of five distinct types. Here we
prove that in the classical framework, with an immovable media of
particles, also the two-dimensional case is more rich: in certain
cases of input of data (height $H$ and radius $r$ of the body) the
problem is well-posed (admit a local minima) without imposing the
restriction $\dot{y}(x) \geq 0$ on the admissible curves
$y(\cdot)$. This is different from the three-dimensional classical
problem or the problem in higher-dimensions, where the restriction
$\dot{y}(x) \geq 0$ is always necessary for the problem to make
sense: without it there exists no strong and no weak local minimum
for Newton's problem of minimal resistance (see \textrm{e.g.}
\cite{CD:Veubeke:1966,CD:CJSilvaMScThesis}). We show that for $H >
\frac{\sqrt 3}{3}r$ the function $\hat{y}(x)= \frac{H}{r}x$ is a
local minimum for the unrestricted Newton's problem of minimal
resistance in dimension two. In the restricted case, while in
dimension three (or higher-orders) the problem has always a unique
solution, we prove that infinitely many different minimizers
appear in dimension two for $r \ge H$. These simple facts seem to
be new in the literature, and never noticed before.

%%%%%%%%%%%%%%%%%%%%%%%%%%%%%%%%%%%%%%%%%

\section{Restricted and unrestricted problems}

In the classical three dimensional Newton's problem of minimal
aerodynamical resistance, the resistance force is given by
$R\left[ \dot{y}(\cdot)\right] = \int_0^r \frac{x}{1 + \dot{y}(x)^2} \, dx$.
Minimization of this functional is a typical problem
of the calculus of variations. Most part
of the old literature wrongly assume the
classical Newton's problem to be ``one of the first applications
of the calculus of variations''. The truth, as Legendre first
noticed in 1788 (see \cite{CD:belloniKawohl:1997}), is that
some restrictions on the derivatives of admissible
trajectories must be imposed: $\dot{y}(x) \geq 0$, $x \in [0, r]$.
The restriction is crucial, because without it there exists no solution,
and the problem suffers from Perron's paradox \cite[\S 10]{CD:MR41:4337}:
since the \emph{a priori} assumption that a
solution exists is not fulfilled, does not make any sense
to try to find it by applying necessary optimality conditions.
It turns out that, with the necessary restriction, the problem is better considered as
an optimal control one (see \cite[p.~67]{Tikhomirov2}
and \cite{CD:arXiv:math.OC/0404237}).
Correct formulation of Newton's problem of minimal resistance
in dimension three is (\textrm{cf. e.g.}
\cite{CD:Veubeke:1966,CD:Tikhomirov:1990}):
\begin{equation*}
\begin{gathered}
\mathcal{R}\left[u(\cdot)\right] = \int_0^r \frac{x}{1 +
u(x)^2} dx \longrightarrow \min \, , \\
\dot{y}(x) = u(x) \, , \quad u(x) \geq 0 \, ,\\
y(0) = 0 \, , \quad y(r) = H \, , \quad H > 0 \, ,
\end{gathered}
\end{equation*}
where we minimize the resistance $\mathcal{R}$ in the class of continuous
functions $y : [0,r] \rightarrow \mathbb{R}$
with piecewise continuous derivative.
Here we consider Newton's problem of minimal resistance in
dimension two (see \cite{CD:arXiv:math.OC/0404237}):
\begin{equation}
\label{eq:R-2-dim}
\begin{gathered}
R\left[ u(\cdot)\right] = \int_0^r \frac{1}{1 + u(x)^2} dx
\longrightarrow \min \, , \\
\dot{y}(x) = u(x) \, , \quad u(x) \in \Omega \, ,\\
y(0) = 0 \, , \quad y(r) = H \, , \quad H > 0 \, .
\end{gathered}
\end{equation}
We consider two cases: (i) unrestricted problem, where no
restriction on the admissible trajectories $y(\cdot)$ other than
the boundary conditions $y(0) = 0$, $y(r) = H$ is considered
($\Omega = \mathbb{R}$); (ii) restricted problem, where the
admissible functions must satisfy the restriction $\dot{y}(x) \geq
0$, $x \in [0, r]$ ($\Omega = \mathbb{R}_0^+$). While for the
classical three-dimensional problem only the restricted problem
admits a minimizer, we prove in \S\ref{sec:UP} that the
two-dimensional problem \eqref{eq:R-2-dim} is more rich: the
unrestricted case also admits a local minimizer when the given
height $H$ of the body is big enough. In \S\ref{sec:RP} we study
the restricted problem. Also in the restricted case the
two-dimensional problem is more interesting: if $r \ge H$, then
infinitely many different minimizers are possible, while in the
classical three-dimensional problem the minimizer is always
unique.

%%%%%%%%%%%%%%%%%%%%%%%%%%%%%%%%%%%

\section{General results for both problems}

The central result of optimal control theory is the Pontryagin
Maximum Principle \cite{CD:MR29:3316b}, which gives a
generalization of the classical necessary optimality
conditions of the calculus of variations.
The following results are valid for both restricted and unrestricted
problems: respectively $\Omega = \mathbb{R}_0^+$ and $\Omega = \mathbb{R}$
in \eqref{eq:R-2-dim}.

\begin{theorem}[Pontryagin Maximum Principle for \eqref{eq:R-2-dim}]
\label{PontMP-d-dim} If $(y(\cdot),u(\cdot))$ is a minimizer of
problem \eqref{eq:R-2-dim}, then there exists a non-zero
pair $(\psi_0,\psi(\cdot))$, where $\psi_0 \leq 0$ is a constant
and $\psi(\cdot) \in PC^1\left([0, r];\,\mathbb{R}\right)$, such
that the following conditions are satisfied for almost all $x$ in
$[0,r]$:
\begin{itemize}
\item [(i)] the Hamiltonian system

\begin{equation*}
\begin{cases}
\dot{y}(x)&=\frac{ \partial{{\cal H}} }{\partial {\psi}}(u(x),
\psi_0, \psi(x)) \quad \text{(control equation $\dot{y} = u$)} \, , \\

\dot{\psi}(x) &= -\frac{ \partial{{\cal H}} }{\partial {y}}(u(x),
\psi_0, \psi(x)) \quad \text{(adjoint system $\dot{\psi} = 0$)} \, ;
\end{cases}
\end{equation*}
\item [(ii)] the maximality condition
\begin{equation}
\label{condMax} {\cal H}(u(x),\psi_0,\psi(x))
={ \mathop {\max}\limits_{u \in \Omega}} {\cal H}(u, \psi_0, \psi(x)) \, ;
\end{equation}
\end{itemize}
where the Hamiltonian ${\cal H}$ is defined by
\begin{equation}
\label{eq:def:ham}
{\cal H}(u,\psi_0,\psi) = \psi_0 \frac{1}{1+u^2}+\psi u \, .
\end{equation}
\end{theorem}
The adjoint system asserts that $\psi(x) \equiv c$,
with $c$ a constant. From the maximality condition it follows that
$\psi_0 \ne 0$ (there are no abnormal extremals for problem \eqref{eq:R-2-dim}).
\begin{proposition}
\label{pont-ext-d-dim}
All the Pontryagin extremals
$\left(y(\cdot), u(\cdot), \psi_0, \psi(\cdot) \right)$ of problem
\eqref{eq:R-2-dim} are normal extremals ($\psi_0 \ne 0$),
with $\psi(\cdot)$ a negative constant: $\psi(x) \equiv
-\lambda$, $\lambda >0$, $x \in [0,r]$.
\end{proposition}
\begin{proof}
The Hamiltonian ${\cal H}$ for problem \eqref{eq:R-2-dim},
$\mathcal{H}\left(u,\psi_0,\psi\right) = \psi_0 \frac{1}{1+u^2} + \psi u$,
does not depend on $y$. Therefore, by
the adjoint system we conclude that
$$\dot{\psi}(x)= - \frac{\partial \mathcal{H}}{\partial y}
\left(u(x),\psi_0,\psi(x) \right)=0\, , $$ that is, $\psi(x)
\equiv c$, $c$ a constant, for all $x \in [0,r]$.
If $c=0$, then $\psi_0 < 0$ (because one can not have both $\psi_0$ and $\psi$ zero)
and the maximality condition \eqref{condMax} simplifies to
\begin{equation}
\label{eq:maxCondAppNP2}
\frac{\psi_0}{1 + u^2(x)}
= { \mathop {\max }\limits_{u \in \Omega}}\left\{\frac{\psi_0}{1 + u^2} \right\} \, .
\end{equation}
From \eqref{eq:maxCondAppNP2} we conclude that the maximum is not achieved
($u \rightarrow \infty$). Therefore $c \neq 0$. Similarly, for $c>0$ the maximum
\begin{equation*}
\frac{\psi_0}{1+u^2(x)} + cu(x)
= {\mathop {\max}\limits_{u \in \Omega}}\left\{\frac{\psi_0}{1+u^2}+ cu\right\}
\end{equation*}
does not exist, and we conclude that $c<0$.
It remains to prove that $\psi_0 \neq 0$. Let us assume
$\psi_0 = 0$. Then the maximality condition reads
\begin{equation}
\label{eq:td2c}
c u(x) = { \mathop {\max }\limits_{u \in \Omega }} \{c u \} \, ,
\quad c < 0 \, .
\end{equation}
For $\Omega = \mathbb{R}$ the maximum does not exist,
and we conclude $\psi_0 \neq 0$.
For $\Omega = \mathbb{R}_0^+$ \eqref{eq:td2c} imply
$u(x) \equiv 0$ and $y(x) \equiv w$, $w$ a
constant ($\dot{y}(x) = u(x)$). This is not possible, given the
boundary conditions $y(0)=0$ and $y(r)=H$ with $H>0$.
Therefore $\psi_0 \neq 0$: there exists no abnormal Pontryagin extremals.
\end{proof}

\begin{remark}
\label{remark-d-dim}
If $\left(y(\cdot), u(\cdot), \psi_0,\psi(\cdot) \right)$
is an extremal, then
$\left(y(\cdot), u(\cdot), \gamma \psi_0, \gamma \psi(\cdot)
\right)$ is also a Pontryagin extremal, for all $\gamma > 0$.
Therefore one can fix, without loss of generality, $\psi_0 = -1$.
\end{remark}
From Proposition~\ref{pont-ext-d-dim} and Remark~\ref{remark-d-dim}
it follows that the Hamiltonian \eqref{eq:def:ham} takes the form
\begin{equation}
\label{eq:HamCNMC}
{\mathcal H}\left(u\right)
= -\frac{1}{1+u^2} - \lambda u \, , \quad \lambda > 0 \, .
\end{equation}

It is not easy to prove the existence of a solution for problem
\eqref{eq:R-2-dim} with classical arguments.
We will use a different approach. We will show,
following \cite{CD:arXiv:math.OC/0404237},
that for problem \eqref{eq:R-2-dim}
the Pontryagin extremals are absolute
minimizers. This means that to solve problem
\eqref{eq:R-2-dim} it is enough to identify its
Pontryagin extremals.

\begin{theorem}
\label{TeormExtPont-d-dim}
Pontryagin extremals for problem \eqref{eq:R-2-dim} are
absolute minimizers.
\end{theorem}
\begin{proof}
Let $\hat{u}(\cdot)$ be a Pontryagin extremal control for problem
\eqref{eq:R-2-dim}. We want to prove that
\begin{equation*}
\int_0^r \frac{1}{1+u^2(x)} dx \ge
\int_0^r \frac{1}{1+\hat{u}^2(x)} dx
\end{equation*}
for any admissible control $u(\cdot)$. Given \eqref{eq:HamCNMC},
we conclude from the maximality condition \eqref{condMax} that
\begin{equation}
\label{cdmaxPN} - \frac{1}{1+\hat{u}^2(x)} - \lambda \hat{u}(x)
\geq - \frac{1}{1+u^2(x)} - \lambda u(x)
\end{equation}
for all $u(\cdot) \in PC \left([0,r], \Omega \right)$.
Having in mind that all the admissible processes $\left(y(\cdot),
u(\cdot)\right)$ of \eqref{eq:R-2-dim} satisfy
\begin{equation*}
\int_0^r u(x) dx = \int_0^r \dot{y}(x) dx = y(r) - y(0) = H \, ,
\end{equation*}
we only need to integrate \eqref{cdmaxPN} to
conclude that $\hat{u}(\cdot)$ is an absolute control minimizer:
\begin{equation*}
\begin{split}
& \int_0^r  \left( - \frac{1}{1+\hat{u}^2(x)} - \lambda
\hat{u}(x) \right) dx \geq  \int_0^r  \left(- \frac{1}{1+u^2(x)}
- \lambda u(x) \right) dx\\
 & \Leftrightarrow \int_0^r  \frac{1}{1+\hat{u}^2(x)} dx
 + \lambda \int_0^r \hat{u}(x) dx \leq \int_0^r \frac{1}{1+u^2(x)} dx
+ \lambda \int_0^r u(x) dx \\
&\Leftrightarrow \int_0^r \frac{1}{1+\hat{u}^2(x)} dx + \lambda H
\leq  \int_0^r  \frac{1}{1+u^2(x)} dx + \lambda H \\
&\Leftrightarrow \int_0^r \frac{1}{1+\hat{u}^2(x)} dx
\leq \int_0^r \frac{1}{1+u(x)^2} dx \, .
\end{split}
\end{equation*}
\end{proof}

Roughly speaking, Theorem~\ref{TeormExtPont-d-dim} reduces the infinite dimension
optimization problem \eqref{eq:R-2-dim} to the study
of a one-dimension maximization problem:
\begin{equation}
\label{eq:opt:dim1}
\max_{u \in \Omega} {\mathcal H}\left(u\right)
= \max_{u \in \Omega} \left\{ -\frac{1}{1+u^2} - \lambda u \right\}
\, , \quad \lambda > 0 \, .
\end{equation}

%%%%%%%%%%%%%%%%%%%%%%%%%%%%%%%%%%%

\section{Unrestricted problem}
\label{sec:UP}

The following standard result of calculus (see \textrm{e.g.} \cite{fenske})
will be used in the sequel.

\begin{theorem}
\label{th:cns:calc:ho}
Let $n \ge 2$ and $\Omega \subseteq \mathbb{R}$ be an open set.
If $f : \Omega \rightarrow \mathbb{R}$ is $n-1$ times differentiable
on $\Omega$ and $n$ times differentiable at some point $a \in \Omega$
where $f^{(k)}(a) = 0$ for $k = 0,\ldots,n-1$ and $f^{(n)}(a) \ne 0$, then:
\begin{itemize}
\item either $n$ is even, and $f(\cdot)$ has an extremum at $a$,
that is a maximum in case $f^{(n)}(a) < 0$ and a minimum in case
$f^{(n)}(a) > 0$;
\item or $n$ is odd, and $f(\cdot)$ does not attain a local extremum at $a$.
\end{itemize}
\end{theorem}

We are considering now the unrestricted two-dimensional
Newton's problem of minimal resistance, that is,
$\Omega = \mathbb{R}$ in \eqref{eq:R-2-dim}.
A necessary (sufficient) condition for $u$
to be a local maximizer for problem \eqref{eq:opt:dim1} is given by
${\mathcal H}'\left(u\right) = 0$ and ${\mathcal H}''\left(u\right) \le 0$
(${\mathcal H}''\left(u\right) < 0$), where
\begin{gather*}
{\mathcal H}'\left(u\right) = {\frac {2 u}{ \left( 1+{u}^{2} \right) ^{2}}}-\lambda \, ,\\
{\mathcal H}''\left(u\right) = -2\,{\frac {3\,{u}^{2}-1}{ \left( 1+{u}^{2} \right) ^{3}}} \, .
\end{gather*}
From the first order condition (maximality condition \eqref{condMax}) it follows that
\begin{equation}
\label{LCN-d-2-naoparam}
\frac{u(x)}{\left(1+u^2(x) \right)^2}=\frac{\lambda}{2}
\Leftrightarrow
\frac{\dot{y}(x)}{\left(1+ \dot{y}^2(x) \right)^2}= \frac{\lambda}{2} \, .
\end{equation}
Using the boundary conditions $y(0)=0$ and $y(r)=H$,
we conclude that $y(x)=\frac{H}{r}x$ ($u = \frac{H}{r}$)
is a local candidate for the solution of the unrestricted problem
($\lambda = {\frac{2{r}^{3} H}{\left( {r}^{2}+{H}^{2} \right)^{2}}}$).
However, by Theorem~\ref{th:cns:calc:ho},
we conclude that such $u$ is a maximizer only when $H > \frac{\sqrt 3}{3}r$.
For $H < \frac{\sqrt 3}{3}r$ the value $u = \frac{H}{r}$ corresponds
to a local minimizer of ${\mathcal H}\left(u\right)$ since
${\mathcal H}''>0$; for $H = \frac{\sqrt 3}{3}r$
function ${\mathcal H}\left(u\right)$ has neither local maximum nor minimum
since ${\mathcal H}''\left(\frac{\sqrt 3}{3}r\right) = 0$ and
${\mathcal H}'''\left(\frac{\sqrt 3}{3}r\right) = -{\frac {27\sqrt{3}}{16}} \ne 0$.

\begin{theorem}
If $H > \frac{\sqrt 3}{3}r$, then
function $y(x)= \frac{H}{r}x$ is a (local) minimum
for the unrestricted problem \eqref{eq:R-2-dim}.
For $H \le \frac{\sqrt 3}{3}r$ the problem has no solution.
\end{theorem}

\begin{remark}
The unrestricted problem \eqref{eq:R-2-dim} does not admit global
minimum. Take indeed, for large values of the parameter $a$, the
control function
\begin{equation*}
\tilde{u}(x) =
\begin{cases}
a & \text{ if } \quad 0 \le x \le \frac{r}{2} + \frac{H}{2a} \\
-a & \text{ if } \quad \frac{r}{2} + \frac{H}{2a} < x \le r \, .
\end{cases}
\end{equation*}
This gives $R[\tilde{u}(\cdot)] = \frac{r}{1 + a^2}$ which
vanishes as $a \rightarrow + \infty$, showing that no global
solution can exist.
\end{remark}

By the symmetry with respect to the $yy$ axis, a local solution to
the unrestricted two-dimensional Newton's problem of minimal
resistance with $H > \frac{\sqrt 3}{3}r$ is a triangle, with value
for resistance $R$ equal to $\frac{r^3}{r^2 + H^2}$.

%%%%%%%%%%%%%%%%%%%%%%%%%%%%%%%%%%%%%%%%%%%%%%%%%%%%%

\section{Restricted problem}
\label{sec:RP}

We now study problem \eqref{eq:R-2-dim} with $\Omega = \mathbb{R}_0^+$.
In this case the optimal control can take values on the boundary
of the admissible set of control values $\Omega$ ($u = 0$).
If the optimal control $u(\cdot)$ is always taking values
in the interior of $\Omega$, $u(x) > 0$ $\forall$ $x \in [0,r]$,
then the optimal solution must satisfy \eqref{LCN-d-2-naoparam}
and it corresponds to the one found in \S\ref{sec:UP}:
\begin{equation}
\label{eq:t}
u(x) = \frac{H}{r} \, , \quad \forall x \in [0,r] \, ,
\end{equation}
with resistance
\begin{equation}
\label{eq:rt}
R = \frac{r^3}{r^2 + H^2} \, .
\end{equation}
We show next that this is solution of the restricted problem
only for $H \ge r$: for $H \le r$ the minimum value for the resistance is
$R = r - \frac{H}{2}$.

It is clear, from the boundary conditions $y(0) = 0$,
$y(r) = H$, $r > 0$, $H > 0$, that $u(x) = 0$, $\forall$
$x \in [0,r]$, is not a possibility: there must exist at least
one non-empty subinterval of $[0,r]$ for which $u(x) > 0$
(otherwise $y(x)$ would be constant, and it would be not possible
to satisfy simultaneously $y(0)=0$ and $y(r)=H$). The simplest
situations are given by
\begin{equation}
\label{eq:i}
u(x) =
\begin{cases}
0 & \text{ if } \quad 0 \le x \le \xi \, , \\
\frac{H}{r - \xi} & \text{ if }  \quad \xi \le x \le r \, ,
\end{cases}
\end{equation}
or
\begin{equation}
\label{eq:f}
u(x) =
\begin{cases}
\frac{H}{\xi} & \text{ if } \quad 0 \le x \le \xi \, , \\
0 & \text{ if }  \quad \xi \le x \le r \, .
\end{cases}
\end{equation}
We get \eqref{eq:t} from \eqref{eq:i} taking $\xi = 0$;
\eqref{eq:t} from \eqref{eq:f} with $\xi = r$.
For \eqref{eq:i} the resistance is given by
$R(\xi) = \xi + \frac{(r-\xi)^3}{(r-\xi)^2 + H^2}$,
that has a minimum value for $\xi = r - H \ge 0$:
$R(r-H) = r - \frac{H}{2}$,
\begin{equation}
\label{eq:io}
u(x) =
\begin{cases}
0 & \text{ if } \quad 0 \le x \le r-H \, , \\
1 & \text{ if }  \quad r-H \le x \le r \, .
\end{cases}
\end{equation}
For $r = H$ \eqref{eq:io} coincides with \eqref{eq:t}; for $r > H$
\begin{equation*}
\left(r - \frac{H}{2}\right) - \left(\frac{r^3}{r^2 + H^2}\right)
= - \frac{H(r-H)^2}{2(r^2+H^2)} < 0 \, ,
\end{equation*}
and \eqref{eq:io} is better than \eqref{eq:t}. Similarly,
for \eqref{eq:f} the resistance is given by
\begin{equation}
\label{eq:rf}
R(\xi) = \frac{\xi^3}{\xi^2 + H^2} + r - \xi \, ,
\end{equation}
that has minimum value for $\xi = H > 0$:
\begin{equation}
\label{eq:fo}
u(x) =
\begin{cases}
1 & \text{ if } \quad 0 \le x \le H \, , \\
0 & \text{ if }  \quad H \le x \le r \, ,
\end{cases}
\end{equation}
$R(H) = r - \frac{H}{2}$,
which coincides with the value for the resistance obtained with \eqref{eq:io}.
If one compares directly \eqref{eq:rt} with \eqref{eq:rf} one get the conclusion
that \eqref{eq:t} is better than \eqref{eq:f} precisely when $r < H$:
\begin{equation}
\label{eq:difrtrf}
\frac{r^3}{r^2 + H^2} - \left(\frac{\xi^3}{\xi^2 + H^2}+r-\xi\right)
= \frac{\xi H^2 \left(r^2 - r\xi - H^2\right)}{\left[(r-\xi)^2+H^2\right](r^2+H^2)} \, ,
\end{equation}
and since $-H^2 \le r^2 - r\xi - H^2 \le r^2 - H^2$, \eqref{eq:difrtrf} is negative
if $r < H$, that is, for $r < H$ \eqref{eq:t} is better than \eqref{eq:f}.
For $r = H$ \eqref{eq:fo} coincide with \eqref{eq:t}, for $r > H$
\eqref{eq:fo} is better than \eqref{eq:t} and as good as \eqref{eq:io}.

We now show that for $r > H$ it is possible to obtain the resistance value
$r - \frac{H}{2}$ from infinitely many other ways,
but no better (no less value) than this quantity.
Generic situation is given by
\begin{equation}
\label{eq:ug}
u_n(x) =
\begin{cases}
0 & \text{ if } \quad \xi_{2i} \le x \le \xi_{2i + 1} \, ,
\quad i = 0,\ldots,n \, ,\\
\frac{\mu_{i+1}-\mu_i}{\xi_{2i+2}-\xi_{2i+1}} &
\text{ if }  \quad \xi_{2i+1} \le x \le \xi_{2i+2} \, ,
\quad i = 0,\ldots,n-1 \, ,
\end{cases}
\end{equation}
where $n \in \mathbb{N}$, $0 = \xi_0 \le \xi_1 \le \cdots \le \xi_{2n+1} = r$,
$0 = \mu_0 \le \mu_1 \le \cdots \le \mu_{n} = H$. We remark that for the
simplest case $n = 1$ \eqref{eq:ug} simplifies to
\begin{equation*}
u_1(x) =
\begin{cases}
0 & \text{ if } \quad 0 \le x \le \xi_1 \, , \\
\frac{H}{\xi_2 - \xi_1} & \text{ if } \quad \xi_1 \le x \le \xi_2 \, , \\
0 & \text{ if }  \quad \xi_2 \le x \le r \, ,
\end{cases}
\end{equation*}
which covers all the previously considered situations:
for $\xi_1 = 0$, $\xi_2 = r$ we obtain \eqref{eq:t};
for $\xi_2 = r$ \eqref{eq:i}; and for $\xi_1 = 0$ one
obtains \eqref{eq:f}. All Pontryagin control extremals
of the restricted problem are of the form \eqref{eq:ug},
and by Theorem~\ref{TeormExtPont-d-dim}
also the minimizing controls. The resistance force $R_n$ associated
with \eqref{eq:ug} is given by
\begin{multline}
\label{eq:Rnug}
R_n\left(\xi_0,\ldots,\xi_{2n+1},\mu_0,\ldots,\mu_n\right) \\
= \sum_{i=0}^{n} \left(\xi_{2i+1}-\xi_{2i}\right)
+ \sum_{i=0}^{n-1} \frac{\left(\xi_{2i+2}-\xi_{2i+1}\right)^3}{\left(\xi_{2i+2}
-\xi_{2i+1}\right)^2 + \left(\mu_{i+1}-\mu_i\right)^2}  \, .
\end{multline}
It is a simple exercise of calculus to see that function \eqref{eq:Rnug}
has three critical points: two of them not admissible, the third one a minimizer.
The first critical point is defined by $\mu_i = 0$, $i = 0,\ldots,n$, which is not
admissible given the fact that $\mu_n = H > 0$. The second critical point is given by
$\mu_i - \mu_{i-1} = \xi_{2i-1} - \xi_{2i}$, $i = 1,\ldots,n$, which is not admissible
since $\mu_i - \mu_{i-1} \ge 0$, $\xi_{2i-1} - \xi_{2i} \le 0$, and
$\mu_i = \mu_{i-1}$, $i = 1,\ldots,n$, is not a possibility given $\mu_n = H > \mu_0 = 0$.
The third critical point is
\begin{equation}
\label{eq:minug}
\mu_i - \mu_{i-1} = \xi_{2i} - \xi_{2i-1}\, , \quad i = 1,\ldots,n \, ,
\end{equation}
which is a minimizer for $H \le r$. Thus, all the minimizing controls
for the restricted two-dimensional problem
with $H \le r$ are of the following form:
\begin{equation}
\label{eq:minContug}
u_n(x) =
\begin{cases}
0 & \text{ if } \quad \xi_{2i} \le x \le \xi_{2i + 1} \, ,
\quad i = 0,\ldots,n \, ,\\
1 & \text{ if }  \quad \xi_{2i+1} \le x \le \xi_{2i+2} \, ,
\quad i = 0,\ldots,n-1 \, ,
\end{cases}
\end{equation}
$n = 1,2,\ldots$, $0 = \xi_0 \le \xi_1 \le \cdots \le \xi_{2n+1} = r$.
For $u_n(x)$ given by \eqref{eq:minContug} the resistance \eqref{eq:Rnug}
reduces to $R_n = r - \frac{H}{2}$, $\forall$ $n \in \mathbb{N}$.

\begin{theorem}
The restricted two-dimensional  Newton's problem of minimal resistance admit always a solution:
\begin{itemize}
\item the unique solution associated to control \eqref{eq:t}, when $H > r$;
\item infinitely many solutions associated to the controls \eqref{eq:minContug}, when $H \le r$.
\end{itemize}
In the case $H > r$ the minimum value for the resistance is $\frac{r^3}{r^2+H^2}$,
otherwise $r - \frac{H}{2}$.
\end{theorem}

%%%%%%%%%%%%%%%%%%%%%%%%%%%%%%%%%%%%%%%%%%%%%%%%%

\section{Conclusion}

Newton's classical problem of minimal resistance offer two
interesting situations to be studied: the problem in dimension
two; and the problem in dimension $d$, $d$ a real number greater
or equal than three. While second situation is well studied in the
literature, and well understood, the first one has been ignored.
In the classical three-dimensional Newton's problem of minimal
resistance, only the problem with restriction $u(x) = \dot{y}(x)
\geq 0$ makes sense (without the restriction the problem has no
local minimum). In the two-dimensional case, we have proved that
the unrestricted case is also a well defined problem when $H >
\frac{\sqrt{3}}{3} r$, the minimum value for the resistance being
$\frac{r^3}{r^2+H^2}$. The local minimizer is a triangle. The
two-dimensional problem with restriction $u(x) = \dot{y}(x) \geq
0$ has always a solution: a unique solution (a triangle) when $H >
r$, with value for resistance equal to the unrestricted case;
infinitely many alternative solutions for $r \ge H$, the minimal
aerodynamical resistance being $r - \frac{H}{2}$.

%%%%%%%%%%%%%%%%%%%%%%%%%%%%%%%%%%%%%%%%%

\section*{Acknowledgments}

This study was proposed to the authors by Alexander Plakhov.
The authors are grateful to him for stimulating discussions.
CS acknowledges the support of the Department
of Mathematics of the University of Aveiro for participation
in the 4th Junior European Meeting; DT the support from
project FCT/FEDER POCTI/MAT/41683/2001.

%%%%%%%%%%%%%%%%%%%%%%%%%%%%%%%%%%%%%%%%%


\begin{thebibliography}{99}

\bibitem{CD:Amaral:1913} I. M. Azevedo do Amaral.
Note sur la solution finie d'un probl\`{e}me de Newton,
Ann. Ac. Pol. Porto, Vol.~8, pp.~207--209, 1913.

\bibitem{CD:belloniKawohl:1997} M. Belloni, B. Kawohl.
A paper of Legendre revisited, Forum Mathematicum,
Vol.~9, pp.~655--668, 1997.

\bibitem{MR56:4953} A.~E. Bryson, Yu~Chi Ho.
Applied optimal control,
Hemisphere Publishing Corp. Washington, D. C.
Optimization, estimation, and control, Revised printing, 1975.

\bibitem{CD:BK:1993} G. Buttazzo, B. Kawohl.
On Newton's problem of minimal resistance,
Math. Intelligencer, Vol.~15, No.~4, pp.~7--12, 1993.

\bibitem{CD:CL:2001} M. Comte, T.~Lachand-Robert.
Newton's problem of the body of minimal resistance
under a single-impact assumption,
Calc. Var. Partial Differential Equations, Vol.~12,
No.~2, pp.~173--211, 2001.

\bibitem{fenske} C. C. Fenske.
Extrema in case of several variables,
The Mathematical Intelligencer, Vol.~25,
No.~1, pp.~49--51, 2003.

\bibitem{CD:Veubeke:1966} B. Fraeijs de Veubeke.
Le probl\`eme de Newton du solide de r\'evolution
pr\'esentant une tra\^\i n\'ee minimum,
Acad. Roy. Belg. Bull. Cl. Sci., Vol.~52,
No.~5, pp.~171--182, 1966.

\bibitem{CD:friction:2002} D. Horstmann, B. Kawohl, P. Villaggio.
Newton's aerodynamic problem in the presence of friction,
NoDEA Nonlinear Differential Equations Appl., Vol.~9, No.~3,
pp.~295--307, 2002.

\bibitem{CD:Kneser:1913} A.~Kneser, E.~Zermelo, H.~Hahn, and M.~Lecat.
Probl\`{e}me de Newton et questions analogues --
Surfaces propulsives, Encyclop\'{e}die des sciences math\'{e}matiques pures
et appliqu\'{e}es, \'{E}dition Française,
Tome II, Vol.~6, Fasc.~1, Calcul des variations,
Paris: Gauthier Villars, Leipzig: B. G. Teubner,
pp.~243--250, 1913.

\bibitem{CD:P1:2003} A.~Yu. Plakhov.
Newton's problem of the body of minimal aerodynamic resistance,
Doklady of the Russian Academy of Sciences, Vol.~390,
No.~3, pp.~1--4, 2003.

\bibitem{CD:arXiv:math.OC/0404194} A. Yu. Plakhov, D. F. M. Torres.
Two-dimensional problems of minimal resistance
in a medium of positive temperature,
Proceedings of the 6th Portuguese Conference
on Automatic Control - Controlo 2004, pp. 488--493, 2004.

\bibitem{math.OC/0407406} A. Yu. Plakhov, D. F. M. Torres.
\textit{Newton's aerodynamic problem in media of chaotically
moving particles}, Sbornik: Mathematics, Volume 196 (2005),
Number 6, pp.~885--933.

\bibitem{CD:MR29:3316b} L.~S.~Pontryagin,
V.~G. Boltyanskii, R.~V. Gamkrelidze, E.~F. Mishchenko.
The mathematical theory of optimal processes, Interscience
Publishers John Wiley \& Sons, Inc.\, New York-London, 1962.

\bibitem{CD:CJSilvaMScThesis} C. J. Silva. Abordagens do C\'{a}lculo
das Variaç\~{o}es e Controlo \'{O}ptimo ao Problema de Newton de
Resist\^{e}ncia M\'{\i}nima, M.Sc. thesis (supervisor:
Delfim F. M. Torres), Univ. of Aveiro, Portugal, June 2005.

\bibitem{CD:Tikhomirov:1990} V. M. Tikhomirov.
Stories about maxima and minima, American
Mathematical Society, Providence, RI, 1990.

\bibitem{Tikhomirov2} V. M. Tikhomirov. Extremal problems -
past and present. In \emph{The Teaching of Mathematics},
Vol.~2, pp.~59-69, 2002.

\bibitem{CD:arXiv:math.OC/0404237} D. F. M. Torres, A. Yu. Plakhov.
Optimal control of Newton-type problems of minimal resistance,
Rend. Semin. Mat. Univ. Politec. Torino Vol.~64 (2006),
no. 1, pp.~79--95.

\bibitem{CD:MR41:4337} L. C. Young.
Lectures on the calculus of variations and optimal control theory,
W. B. Saunders Co., Philadelphia, 1969.

\end{thebibliography}
\end{document}